\documentclass[12pt,letterpaper,final,twoside,leqno]{amsart}
\usepackage{amsmath,amssymb,amsthm,amsfonts,amscd,amsopn}
\usepackage{hyperref,datetime}
\usepackage[usenames]{color}
\usepackage{eucal, mathrsfs}
\usepackage{xspace}
\usepackage[all]{xy}\xyoption{dvips}
\usepackage{comment} 
\usepackage{enumitem}
\usepackage{times}
\usepackage[T1]{fontenc}
\newcommand{\skpagesize}{
\textwidth= 6.25in
\textheight=8.75in
\voffset-.5in
\hoffset-.5in
\marginparwidth=56pt
}
\definecolor{pinktwo}{RGB}{94,52,121}
\newcommand\sandor[1]{{\color{pinktwo}{\tiny\sf S\'andor\,\mdyydate{#1}}}}
\newcommand\sideremark[1]{%
\normalmarginpar
\marginpar
[
\hskip .45in
\begin{minipage}{.75in}
\tiny #1
\end{minipage}
]
{
\hskip -.075in
\begin{minipage}{.75in}
\tiny #1
\end{minipage}
}}
\newcommand\rsideremark[1]{
\reversemarginpar
\marginpar
[
\hskip .45in
\begin{minipage}{.75in}
\tiny #1
\end{minipage}
]
{
\hskip -.075in
\begin{minipage}{.75in}
\tiny #1
\end{minipage}
}}
\newcommand\silentremark[1]{}
\newtheoremstyle{bozont}{3pt}{3pt}%
     {\itshape}%         Body font
     {}%         Indent amount (empty = no indent, \parindent = para indent)
     {\bfseries}% Thm head font
     {.}%        Punctuation after thm head
     {.5em}%     Space after thm head (\newline = linebreak)
     {\thmname{#1}\thmnumber{ #2}\thmnote{ \rm #3}}%         Thm head spec
\newtheoremstyle{bozont-sf}{3pt}{3pt}%
     {\itshape}%         Body font
     {}%         Indent amount (empty = no indent, \parindent = para indent)
     {\sffamily}% Thm head font
     {.}%        Punctuation after thm head
     {.5em}%     Space after thm head (\newline = linebreak)
     {\thmname{#1}\thmnumber{ #2}\thmnote{ \rm #3}}%         Thm head spec
\newtheoremstyle{bozont-sc}{3pt}{3pt}%
     {\itshape}%         Body font
     {}%         Indent amount (empty = no indent, \parindent = para indent)
     {\scshape}% Thm head font
     {.}%        Punctuation after thm head
     {.5em}%     Space after thm head (\newline = linebreak)
     {\thmname{#1}\thmnumber{ #2}\thmnote{ \rm #3}}%         Thm head spec
\newtheoremstyle{bozont-remark}{3pt}{3pt}%
     {}%         Body font
     {}%         Indent amount (empty = no indent, \parindent = para indent)
     {\scshape}% Thm head font
     {.}%        Punctuation after thm head
     {.5em}%     Space after thm head (\newline = linebreak)
     {\thmname{#1}\thmnumber{ #2}\thmnote{ \rm #3}}%         Thm head spec
\newtheoremstyle{bozont-def}{3pt}{3pt}%
     {}%         Body font
     {}%         Indent amount (empty = no indent, \parindent = para indent)
     {\bfseries}% Thm head font
     {.}%        Punctuation after thm head
     {.5em}%     Space after thm head (\newline = linebreak)
     {\thmname{#1}\thmnumber{ #2}\thmnote{ \rm #3}}%         Thm head spec
\newtheoremstyle{bozont-reverse}{3pt}{3pt}%
     {\itshape}%         Body font
     {}%         Indent amount (empty = no indent, \parindent = para indent)
     {\bfseries}% Thm head font
     {.}%        Punctuation after thm head
     {.5em}%     Space after thm head (\newline = linebreak)
     {\thmnumber{#2.}\thmname{ #1}\thmnote{ \rm #3}}%         Thm head spec
\newtheoremstyle{bozont-reverse-sc}{3pt}{3pt}%
     {\itshape}%         Body font
     {}%         Indent amount (empty = no indent, \parindent = para indent)
     {\scshape}% Thm head font
     {.}%        Punctuation after thm head
     {.5em}%     Space after thm head (\newline = linebreak)
     {\thmnumber{#2.}\thmname{ #1}\thmnote{ \rm #3}}%         Thm head spec
\newtheoremstyle{bozont-reverse-sf}{3pt}{3pt}%
     {\itshape}%         Body font
     {}%         Indent amount (empty = no indent, \parindent = para indent)
     {\sffamily}% Thm head font
     {.}%        Punctuation after thm head
     {.5em}%     Space after thm head (\newline = linebreak)
     {\thmnumber{#2.}\thmname{ #1}\thmnote{ \rm #3}}%         Thm head spec
\newtheoremstyle{bozont-remark-reverse}{3pt}{3pt}%
     {}%         Body font
     {}%         Indent amount (empty = no indent, \parindent = para indent)
     {\sc}% Thm head font
     {.}%        Punctuation after thm head
     {.5em}%     Space after thm head (\newline = linebreak)
     {\thmnumber{#2.}\thmname{ #1}\thmnote{ \rm #3}}%         Thm head spec
\newtheoremstyle{bozont-def-reverse}{3pt}{3pt}%
     {}%         Body font
     {}%         Indent amount (empty = no indent, \parindent = para indent)
     {\bfseries}% Thm head font
     {.}%        Punctuation after thm head
     {.5em}%     Space after thm head (\newline = linebreak)
     {\thmnumber{#2.}\thmname{ #1}\thmnote{ \rm #3}}%         Thm head spec
\newtheoremstyle{bozont-def-newnum-reverse}{3pt}{3pt}%
     {}%         Body font
     {}%         Indent amount (empty = no indent, \parindent = para indent)
     {\bfseries}% Thm head font
     {}%        Punctuation after thm head
     {.5em}%     Space after thm head (\newline = linebreak)
     {\thmnumber{#2.}\thmname{ #1}\thmnote{ \rm #3}}%         Thm head spec
\newtheoremstyle{bozont-def-newnum-reverse-plain}{3pt}{3pt}%
   {}% Body font
   {}% Indent amount (empty = no indent, \parindent = para indent)
   {}% Thm head font
   {}% Punctuation after thm head
   {.5em}% Space after thm head (\newline = linebreak)
   {\thmnumber{\!(#2)}\thmname{ #1}\thmnote{ \rm #3}}% Thm head spec
\newcounter{thisthm} 
\newcounter{thissection} 
\newcommand{\ilabel}[1]{%
  \newcounter{#1}%
  \setcounter{thissection}{\value{section}}%
  \setcounter{thisthm}{\value{thm}}%
  \label{#1}}
\newcommand{\iref}[1]{%
  (\the\value{thissection}.\the\value{thisthm}.\ref{#1})}
\numberwithin{figure}{section}

\newlength{\swidth}
\skpagesize
\renewcommand\sideremark[1]{}\renewcommand\rsideremark[1]{}
\pdfpagewidth=\paperwidth \pdfpageheight=\paperheight
\providecommand{\bysame}{\leavevmode\hbox to3em{\hrulefill}\thinspace}
\usepackage{hyperref}
\usepackage[latin1]{inputenc}
\usepackage[all]{xy}
\usepackage{indentfirst}
\usepackage{fancyhdr}
\usepackage{textcomp}
\usepackage{graphicx}
\usepackage{verbatim}
\usepackage{stmaryrd}
\usepackage{pifont}
\usepackage{eucal}
\usepackage{tikz}
\usepackage{mathrsfs}

\theoremstyle{bozont}
\newtheorem{thm}{Theorem}[section]
\newtheorem{theorem}[thm]{Theorem}
\newtheorem{lemma}[thm]{Lemma}
\newtheorem{proposition}[thm]{Proposition}
\newtheorem{corollary}[thm]{Corollary}
\newtheorem{propdef}[thm]{Proposition-Definition}

\theoremstyle{bozont-def}
\newtheorem{definition}[thm]{Definition}
\newtheorem{remark}[thm]{Remark}
\newtheorem{example}[thm]{Example}

\newtheorem{question}[thm]{Question}

\newtheorem{convention}[thm]{Convention}
\newtheorem{de}[thm]{Definition}

\newtheorem{no}[thm]{Notation}

\newtheorem{qsn}[thm]{Question}

\newcommand\eps{\varepsilon}
\renewcommand\ge{\geqslant}  % schr�ge Variante
\renewcommand\geq{\geqslant}  % schr�ge Variante
\renewcommand\leq{\leqslant}  % schr�ge Variante

\newcommand\be{\begin{eqnarray*}}
\newcommand\ee{\end{eqnarray*}}

\newcommand\Q{\mathbb Q}

\newcommand\R{\mathbb R}
\renewcommand\P{\mathbb P}
\newcommand\sO{\mathscr{O}}

\newcommand\with{\,\,\vrule\,\,}
\newcommand\set[1]{\left\{#1\right\}}
\newcommand\newop[2]{\def#1{\mathop{\rm #2}\nolimits}}
\newop\log{log}
\newop\ord{ord}
\newop\Gal{Gal}
\newop\SL{SL}
\newop\mult{mult}
\newop\Eff{Eff}
\newop\Mov{Mov}
\newop\Amp{Amp}
\newop\Null{Null}
\newop\Supp{Supp}
\newop\Id{Id}
\newop\rk{rk}
\newop\Sym{Sym}
\newop\mov{mov}
\newop\nonnef{nonnef}

\newcommand{\deq}{\ensuremath{ \stackrel{\textrm{def}}{=}}}
\newcommand{\equ}{\ensuremath{\,=\,}}
\newcommand{\st}[1]{\ensuremath{ \left\{ #1 \right\} }}

\newcommand{\shq}{\mathcal{Q}}
\newcommand{\shb}{\mathcal{B}}

\newcommand{\she}{{E}}
\newcommand{\sha}{{A}}
\newcommand{\shh}{{H}}
\newcommand{\shl}{{L}}

\newcommand\lra{\longrightarrow}
\newcommand\op{\sO_{\P(\she)}(1)}
\newcommand{\HH}[3]{\ensuremath{H^{#1}\left(#2,#3\right)}}

\newcommand{\C}{\mathbb{C}}

\begin{document}

\title{On positivity and base loci of vector bundles} %{Notes on bigness}

\author[Bauer, Kov\'acs, K\"uronya, Mistretta, Szemberg, Urbinati]{Thomas Bauer,
  S\'andor J Kov\'acs, Alex K\"uronya, Ernesto Carlo Mistretta, Tomasz Szemberg, Stefano
  Urbinati}

\address{TB: Fachbereich Mathematik und Informatik, Philipps-Universit\"at Marburg, Hans-Meerwein-Stra{\ss}e, Lahnberge, 35032 Marburg, Germany}
\address{SJK: University of Washington, Department of Mathematics, Seattle, WA 98195, USA}
\address{AK: Budapest University of Technology and Economics, Mathematics Institute, Egry J\'ozsef u. 1, H-1111 Budapest, Hungary}
\address{ECM: Universit\`a degli Studi di Padova, Dipartimento di Matematica, ROOM 619, Via Trieste 63, 35121 Padova, Italy}
\address{TSz: Pedagogical Univarsity Cracow, Podchor\c a\.zych 2, 30-084 Krak\'ow, Poland}
\address{SU: Universit\`a degli Studi di Padova, Dipartimento di Matematica, ROOM 608, Via Trieste 63, 35121 Padova, Italy}

\email{TB: tbauer@mathematik.uni-marburg.de}
\email{SJK: skovacs@uw.edu}
\email{AK: alex.kuronya@math.bme.hu}
\email{ECM: ernesto.mistretta@unipd.it}
\email{TSz: tomasz.szemberg@gmail.com}
\email{SU: urbinati.st@gmail.com}

\date{\usdate\today} %{\it Internal notes of \today}

\maketitle
\thispagestyle{empty}

%*****************************************************************************

\silentremark
{\center{\sandor{060914}: \small\tt\bf PLEASE!! When you add something, try to follow the
  notational conventions already in the paper. I corrected a bunch of
  inconsistencies, and I don't want to do it again. Please look at the rest of the
  paper how things are denoted. E.g., intersection numbers with ``cdot'' but without
  enclosing parantheses. There were three different notations in different parts!}}

\silentremark{\sandor{060914}: I removed ``Viehweg'' from ``weakly positive'' everywhere
  I found, because I don't believe there is any other notion of weak positivity. Same
  for Miyaoka and AEA.}

\section{Introduction}

The aim of this note is to shed some light on the relationships among  some notions of positivity for vector bundles that arose
in recent decades.

Positivity properties of line bundles  have  long played a major role in projective geometry;
they have once again become a center of attention recently, mainly in relation with advances in birational geometry,
especially in the framework of the Minimal Model Program. Positivity of line bundles has often
been studied in conjunction with numerical invariants and various kinds of asymptotic base loci (see
for example \cite{ELMNP} and \cite{BDPP}).

At the same time, many positivity notions have been introduced for vector bundles of higher rank,
generalizing some of the properties that hold for line bundles.  While the situation in rank one is well-understood, at least
as far as the interdepencies between the various positivity concepts is concerned, we are quite far from an analogous state of
affairs for vector bundles in general.

In an attempt to generalize bigness for the higher rank case, some positivity properties have been put forward
by Viehweg (in the study of fibrations in curves, \cite{vie}), and Miyaoka (in the context of surfaces, \cite{Miy83}), and are known to
be different from the generalization given by using the tautological line bundle on
the projectivization of the considered vector bundle (cf. \cite{PAG}).  The
differences between the various definitions of bigness are already present in the
works of Lang concerning the Green-Griffiths conjecture (see \cite{lang}).

Our purpose  is to study  several of the positivity notions studied for vector
bundles with some notions of asymptotic base loci that can be defined on the
variety itself, rather than on the projectivization of the given vector bundle.  We
relate some of the different notions conjectured to be equivalent with the help of
these base loci, and we show that these can help simplify the various relationships
between the positivity properties present in the literature.

In particular, we define augmented and restricted base loci $\mathbb{B}_+ (E)$ and
$\mathbb{B}_- (E)$ of a  vector bundle $E$ on the variety $X$, as generalizations of the corresponding
notions studied extensively for line bundles.  As it turns out, the asymptotic base loci defined here behave well with
respect to the natural map induced by the projectivization of the vector bundle $E$,
as shown in Section~\ref{asvect}.

The relationship between these base loci with the positivity notions appearing in the
literature goes as follows:

\begin{theorem}

  Let $X$ be a smooth projective variety and $E$ a vector bundle on $X$. Then:

\begin{enumerate}

\item $E$ is \emph{ample} if and only if ${\mathbb{B}_+ (E)}  = \emptyset$;

\item $E$ is \emph{nef} if and only if ${\mathbb{B}_- (E)}  = \emptyset$;

\item $E$ is \emph{pseudo-effective} if and only if ${\mathbb{B}_- (E)} \neq X$;
  \ilabel{changed} %
  \silentremark{\sandor{060914}: changed ``$\emptyset$'' to ``$X$'' in \iref{changed}}

\item $E$ is %Viehweg
  \emph{weakly positive} if and only if $\overline{\mathbb{B}_- (E)} \neq X$ (see
  Section~\ref{restricted});

\item $E$ is \emph{V-big} if and only if ${\mathbb{B}_+ (E)} \neq X$ (see section
  \ref{augmented});

\item Assume that $E$ is a nef vector bundle. Then $E$ is %Miyaoka
  \emph{almost everywhere ample} if and only if ${\mathbb{B}_+ (E)} \neq X$ (cf.\
  Section~\ref{badaea}).

\end{enumerate}

\end{theorem}

The paper is organized as follows: in sections \ref{def} and \ref{asvect} we give the
definition and basic properties of the asymptotic base loci for vector bundles, and
relate these loci with the ones on the projectivizations. In Section~\ref{line} we recall the various positivity properties for line bundles
and their relationship with asymptotic base loci. Section~\ref{restricted} is devoted to a  study of  positivity properties of vector bundles
related to the restricted base locus, while   Section~\ref{augmented} is given over to an investigation of connection between
positivity properties of  vector bundles related and  augmented base loci.
In sections \ref{sectionaea} and \ref{badaea} we study almost everywhere ampleness
and relate it to V-bigness.

\section{Definitions and first properties}
\label{def}

\begin{convention}\label{conv}
  Throughout the paper we are working with vector bundles of finite rank, but for
  various reasons we find it  more convenient to work with the associated sheaf of sections
  which is a locally free coherent $\sO_X$-module. We will follow the usual abuse of
  terminology and while exclusively using this associated sheaf, we will still call it a
  vector bundle. If, rarely, we want to refer to a vector bundle and mean a vector
  bundle we will call it the \emph{total space} of the vector bundle.

  We will also work with line bundles, which of course refers to a locally free sheaf
  of rank $1$. For a line bundle $L$ we will denote by $c_1(L)$ the associated Weil
  divisor on $X$.
\end{convention}

With that convention fixed we are making the following notation that we will use
through the entire paper:

\begin{no}
  Let $X$ be a smooth projective variety over the complex numbers, and $E$ a vector
  bundle (i.e., according to \ref{conv} really a locally free sheaf) over $X$. For a
  point $x\in X$, $E_x=E\otimes_{\sO_X}\sO_{X,x}$ denotes the stalk of $E$ at the
  point $x$ and $E(x)=E\otimes_{\sO_X} \kappa(x)$ where $\kappa(x)$ is the residue
  field at $x$. Clearly, $E(x)$ is the fiber of the total space of $E$ over the point
  $x$. In particular, $E(x)$ is a vector space of dimension $r=\rk E$.
\end{no}

\begin{definition}

  We define the \emph{base locus} of $E$ (over $X$) as the subset
\[
\mathrm{Bs} (E) := \{ x \in X ~|~ H^0(X, E) \to E(x) ~\textrm{is not surjective} \}\ ,
\]
and the \emph{stable base locus} of $E$ (over $X$) as
\[
\mathbb{B} (E) := \bigcap_{m>0} \mathrm{Bs}(\Sym^m E)\ .
\]
\end{definition}

\begin{remark} The assertions below follow immediately from the definition:
  \begin{enumerate}
  \item As $\mathrm{Bs} (E) = {\rm Bs} (\mathrm{Im} (\bigwedge^{\rk E} H^0(X, E) \to
    H^0(X, \det E) ))$, these loci are closed subsets, and carry a natural scheme
    structure.
  \item For any positive integer $c>0$, $\mathbb{B}(E) = \mathbb{B}(\Sym^c E)$, and
    the same holds for $\mathbb{B}_-$ and $\mathbb{B}_+$. \ilabel{base-of-sym}
  \end{enumerate}
\end{remark}

\begin{remark}
The rank of the natural linear map $H^0(X, E) \to E(x)$ induces a stratification of $X$ into locally closed subsets.
\end{remark}

\begin{definition}

  Let $r = p/q \in \mathbb{Q}_{>0}$ be a positive rational number, and $A$ a line
  bundle on $X$.  We will use the following notation:
  \begin{align*}
    \mathbb{B}(E + rA) &:= \mathbb{B}(\Sym^q E \otimes A^{p}), \quad\text{and}\\
    \mathbb{B}(E - rA) &:= \mathbb{B}(\Sym^q E \otimes A^{-p}) \ .
  \end{align*}

  Note that if $r=p'/q'$ is another representation of $r$ as a fraction, then
  $q'p=p'q$, hence
  $$
  \Sym^{q'}(\Sym^q E \otimes A^p)\simeq \Sym^{q'q} E \otimes A^{q'p}\simeq
  \Sym^{q'q} E \otimes A^{p'q}\simeq \Sym^{q}(\Sym^{q'} E \otimes A^{p'}),
  $$
  and therefore, by \iref{base-of-sym}, $\mathbb{B}(\Sym^q E \otimes
  A^p)=\mathbb{B}(\Sym^{q'} E \otimes A^{p'})$ and hence $\mathbb{B}(E + rA)$ is
  well-defined. A similar argument shows that $\mathbb{B}(E - rA)$ is also
  well-defined.

Let $A$ be an ample line bundle on $X$,
we define the \emph{augmented base locus} of $E$ as
\[
\mathbb{B}_+^A (E):= \bigcap_{r \in \mathbb{Q}^{>0}} \mathbb{B}(E - r A ) \ ,
\]
and the \emph{restricted base locus} of $E$ as
\[
\mathbb{B}_-^A (E):= \bigcup_{r \in \mathbb{Q}^{>0}} \mathbb{B}(E + r A )\ .
\]

\end{definition}

\begin{remark}

The definitions above yield the following properties:

\begin{enumerate}

\item The loci \( \mathbb{B}_+^A (E)\) and \( \mathbb{B}_-^A (E)\) do not depend on
  the choice of the ample line bundle $A$, so we can write
  \( \mathbb{B}_+ (E)\) and \( \mathbb{B}_- (E)\) for the augmented and restricted base locus of $E$, respectively.

\item For any $r_1 > r_2 > 0$ we have \( \mathbb{B} (E + r_1 A) \subseteq \mathbb{B}
  (E + r_2 A) \) and \( \mathbb{B} (E - r_2 A) \subseteq \mathbb{B} (E - r_1 A) \).

\item In particular, for any $\eps >0$ we have \( \mathbb{B} (E + \eps A) \subseteq
  \mathbb{B} (E) \subseteq \mathbb{B} (E - \eps A) \).

\item Therefore we have that
  $$
    \mathbb{B}_+ (E):= \bigcap_{q \in \mathbb{N}}
    \mathbb{B}(E - (1/q) A ) \qquad \text{and} \qquad
    \mathbb{B}_- (E):= \bigcup_{q \in \mathbb{N}} \mathbb{B}(E + (1/q) A ).
  $$

\item It follows that $\mathbb B_+(E)$ is closed, but even for line bundles, the
  locus $\mathbb{B}_- (E)$ is not closed in general: Lesieutre \cite{lesieutre}
  proved that this locus can be a proper dense subset of $X$, or a proper dense
  subset of a divisor of $X$.

\end{enumerate}

\end{remark}

\section{Asymptotic invariants for vector bundles}
\label{asvect}

In the following sections we will relate  augmented and restricted base loci for
vector bundles to various positivity notions found in the literature. In order to
achieve  a better understanding of these positivity properties and the relations
between them, it is necessary to investigate the dependence of
asymptotic base loci for vector bundles,  and the corresponding loci of the tautological quotient line
bundles on the appropriate projectivizations.

Let $E$ be a vector bundle on a smooth projective variety $X$, $\pi \colon
\mathbb{P}(E) \to X$ the projective bundle of rank one quotients of $E$, and
$\sO_{\mathbb{P}(E)}(1)$ the universal quotient of $\pi^* E$ on $\mathbb{P}(E)$.
Then we immediately have
\[
\pi (\mathbb{B} (\sO_{\mathbb{P}(E)}(1))) \subseteq \mathbb{B} (E) \ .
\]
In fact, if the evaluation map $H^0(X, E) \otimes \sO_X \to E$ is surjective over a
point $x \in X$, then the map
\[
H^0(\mathbb P(E), \sO_{\mathbb P(E)}(1)) \otimes\sO_{\mathbb{P}(E)} = H^0(X, E) \otimes \sO_{\mathbb{P}(E)} \to \pi^* E
\twoheadrightarrow \sO_{\mathbb{P}(E)}(1)
\]
is surjective over any point $y \in
\mathbb{P}(E)$ such that $\pi(y) = x$, and a similar argument applies to
$\Sym^m E$.

More precisely, we have $\pi (\mathrm{Bs} (\sO_{\mathbb{P}(E)}(1))) = \mathrm{Bs} (E)$: if a
point $x \in X$ does lie in $\mathrm{Bs} (E)$, then the image of the map $H^0(X, E)
\to E(x)$ is contained in some hyperplane $H \subset E(x)$, where the  hyperplane $H$
corresponds to a point $y \in \pi^{-1} (x)$  contained in $\mathrm{Bs}(\sO_{\mathbb{P}(E)}(1))$.

It is not clear whether the inclusion $\pi(\mathbb{B} (\sO_{\mathbb{P}(E)}(1)))
\subseteq \mathbb{B} (E) $ of  stable loci is strict in general.  However, as we
will show right below, some useful connections  rely on properties
of  augmented and restricted base loci, which exhibit  a more predictable
behavior with respect to the map $\pi$.

\begin{proposition}
\label{bmin}

Let $E$ be a vector bundle on a smooth projective variety $X$,
$\pi \colon \mathbb{P}(E) \to X$ the projective bundle of one dimensional quotients of $E$,
and $\sO_{\mathbb{P}(E)}(1)$ the universal quotient of $\pi^* E$ on $\mathbb{P}(E)$.
Then
\[
\pi(\mathbb{B}_- (\sO_{\mathbb{P}(E)}(1))) =  \mathbb{B}_- (E) ~.
\]

\end{proposition}

\begin{proof}

  Let us fix $H \in \mathrm{Pic} (X)$, a sufficiently ample line bundle such that
  $A:= \sO_{\mathbb{P}(E)}(1) \otimes \pi^* H$ is very ample on $\mathbb{P}(E)$.
  Then
  \[
  \mathbb{B}_- (\sO_{\mathbb{P}(E)}(1)) = \bigcup_{a \in \mathbb{N}}
  \mathbb{B}(\sO_{\mathbb{P}(E)}(a) \otimes A)= \bigcup_{a \in \mathbb{N}}
  \Big(\bigcap_{b\in \mathbb{N}} \mathrm{Bs}(\sO_{\mathbb{P}(E)}(ab) \otimes A^{b})
  \Big)
  \]
  and
  \[
  \mathbb{B}_- (E) = \bigcup_{a\in \mathbb{N}} \Big(\bigcap_{b\in \mathbb{N}}
  \mathrm{Bs} (\Sym^{ab} E \otimes H^{b}) \Big).
  \]

  $(\subseteq)$ The easier inclusion is \( \pi(\mathbb{B}_- (\sO_{\mathbb{P}(E)}(1)))
  \subseteq \mathbb{B}_- (E).  \) In order to show this, suppose that $x \in X$ and
  that $x \notin \mathbb{B}_- (E)$.  Then for any integer $a >0$ there exists a $b>0$
  such that the vector bundle $\Sym^{ab} E \otimes H^b$ is generated by its global
  sections at $x \in X$.  Then for all $a > 0$ the line bundle $\sO(2(a-1)b) \otimes
  A^{2b} = \sO (2ab) \otimes \pi^* H^{2b}$ which is a quotient of $\pi^* (\Sym^{2ab}
  E \otimes H^{2b})$ is generated by its global sections (defined over the whole
  space $\mathbb{P} (E)$) on any point of the fibre $\pi^{-1} (x)$, so the fibre
  $\pi^{-1} (x)$ is contained in the complement of $\mathbb{B}_-
  (\sO_{\mathbb{P}(E)}(1))$.

  $(\supseteq )$ Let us show now that \( \pi(\mathbb{B}_- (\sO_{\mathbb{P}(E)}(1)))
  \supseteq \mathbb{B}_- (E) \): Let $x \in \mathbb{P}(E)$ be a point such that $x
  \notin \pi(\mathbb{B}_- (\sO_{\mathbb{P}(E)}(1))) $.  Then for any $a>0$ there
  exists a $b>0$ such that $\sO_{\mathbb{P}(E)}(2(a-1)b) \otimes A^{b}$ is generated on
  any point $y \in \pi^{-1} (x)$ by its global sections (defined on the whole
  $\mathbb{P}(E)$).

  Then the line bundle
  \begin{align*}
    L := \sO_{\mathbb{P}(E)} (2ab) \otimes \pi^* H^{2b} &\simeq \sO_{\mathbb{P}(E)} (2(a-1)b)
  \otimes \sO_{\mathbb{P}(E)} (b) \otimes \pi^* H^{b} \otimes \sO_{\mathbb{P}(E)} (b)
  \otimes \pi^* H^{b} \simeq\\
  &\simeq \Big (\sO_{\mathbb{P}(E)}(2(a-1)b) \otimes A^{b} \Big ) \otimes A^{b}
  \end{align*}
  is the product of a line bundle which is generated by global sections (on
  $\mathbb{P}(E)$) at any point of the fiber $\mathbb P_x := \pi^{-1} (x) = \mathbb{P}(E(x))$
  with a very ample line bundle, so its global sections (on $\mathbb{P}(E)$) define a
  closed immersion of $\mathbb P_x$ into a projective space.
%  \begin{sloppypar}
    In other words, the linear system $H^0(\mathbb{P}(E), L)$ defines a rational
    map \hbox{ $\varphi \colon \mathbb{P}(E) \dashrightarrow
      \mathbb{P}(H^0(\mathbb{P(E)}, L)) =\mathbb{P}^N$} which is a regular immersion
    on $\mathbb P_x$.  Then in particular, for $m \gg 0$ the multiplication map $\Sym^m
    H^0(\mathbb{P(E)}, L) \to H^0(Y , L|_{Y})$ is surjective.
%  \end{sloppypar}

  It follows that the map $\pi_* (\Sym^m H^0(\mathbb{P}(E), L) \otimes
  \sO_{\mathbb{P}(E)}) \to \pi_*(L^{m})$ is surjective at the point $x \in X$.  As
  $\pi_*(L^{m}) = \Sym^{2abm} E \otimes H^{2bm}$ we may conclude that for any $a>0$
  and $m$ large enough the vector bundle $\Sym^{2abm} E \otimes H^{2bm}$ is generated
  at $x$ by its global sections, hence $x \notin \mathbb{B}_- (E)$.
\end{proof}

The analogous claim  holds for  augmented base locus, with a similar proof.

\begin{proposition}
  \label{bplus}

  Let $E$ be a vector bundle on a smooth projective variety $X$, with the same
  notations as in Proposition \ref{bmin}, we have
  \[
  \pi(\mathbb{B}_+ (\sO_{\mathbb{P}(E)}(1))) =  \mathbb{B}_+ (E) ~.
  \]

\end{proposition}

\begin{proof}

%  As above,
  Let $H \in \mathrm{Pic} (X)$ be a sufficiently ample line bundle
  such that $A:= \sO_{\mathbb{P}(E)}(1) \otimes \pi^* H$ is very ample on
  $\mathbb{P}(E)$.  Then
  \[
  \mathbb{B}_+ (\sO_{\mathbb{P}(E)}(1)) = \bigcap_{a>0}
  \mathbb{B}(\sO_{\mathbb{P}(E)}(a) \otimes A^{-1})= \bigcap_{a>0} \Big(\bigcap_{b>0}
  \mathrm{Bs}(\sO_{\mathbb{P}(E)}(ab) \otimes A^{-b}) \Big)
  \]
  and
  \[
  \mathbb{B}_+ (E) = \bigcap_{a>0} \Big(\bigcap_{b>0} \mathrm{Bs} (\Sym^{ab} E
  \otimes H^{-b}) \Big).
  \]

  In order to show that \( \pi(\mathbb{B}_+ (\sO_{\mathbb{P}(E)}(1))) \subseteq
  \mathbb{B}_+ (E), \) observe that if $\Sym^{ab}E \otimes H^{-b}$ is globally
  generated at a point $x \in X$, then $\pi^* \Sym^{ab}E \otimes \pi^*H^{-b}$ is
  generated at all points in $\pi^{-1}(x)$, hence $\sO_{\mathbb{P}(E)} (ab) \otimes
  \pi^*H^{-b} = \sO_{\mathbb{P}(E)} ((a+1)b) \otimes A^{-b} $ is globally generated
  at all points in $\pi^{-1}(x)$.

  To show the other inclusion, set $U = X \setminus \pi(\mathbb{B}_+
  (\sO_{\mathbb{P}(E)}(1)))$ and observe that $(\sO_{\mathbb{P}(E)} (ab) \otimes
  A^{-b})$ is generated by its global sections at the points of $\pi^{-1} (U)$ for
  $a$ and $b$ sufficiently large.  Let us consider $b>0$ a sufficiently large
  positive integer, $a = (b-1)k >0$ a sufficiently large multiple of $b-1$, and set
  $c:= ((a-1)b +1)/(b-1) = kb-1 = a + k-1$.  Finally, let $L$ be the line bundle
  \begin{multline*}
    L:= \sO_{\mathbb{P}(E)} (c (b-1)) \otimes \pi^*H^{-(b-1)} = \sO_{\mathbb{P}(E)}
    ((a-1)b +1) \otimes \pi^*H^{-(b-1)}\simeq \\
    \simeq \Big( \sO_{\mathbb{P}(E)} (ab) \otimes \big( \sO_{\mathbb{P}(E)} (-1)
    \otimes \pi^*H^{-1} \big)^{b} \Big) \otimes (\sO_{\mathbb{P}(E)} (1) \otimes
    \pi^*H) \simeq (\sO_{\mathbb{P}(E)} (ab) \otimes A^{-b}) \otimes A.
  \end{multline*}

  Now for $b$ and $k$ large enough $L$ is the product of the very ample line bundle
  $A$ with a line bundle which is generated by global sections on $\pi^{-1} (U)$, so
  it is very ample on the open subset $\pi^{-1} (U)$. Furthermore, we have that
  $\pi_* (L) = \Sym^{c(b-1)} E \otimes H^{-(b-1)}$ and so we can apply the same
  argument as in the proof of Proposition~\ref{bmin} to finish the proof.
\end{proof}

\section{Positivity properties for line bundles}
\label{line}

We recall here how  augmented and restricted base loci are involved with various
positivity notions of line bundles, as well as  loci defined by negative curves:

\begin{de}
  Let $L$ be a line bundle on a smooth projective variety $X$.  Fix an ample line
  bundle $A$ and a rational number $\eps>0$. We define
  \begin{enumerate}
  \item ${T_\eps^A=\overline{\set{ x \with \exists\, C\subseteq X \textrm{ curve on $X$
            s.t. } x \in C, c_1(L)\cdot C<\eps\cdot c_1(A)\cdot C}}}$
    to be the  \emph{non-AEA locus of $L$ with respect to $A$ and $\eps$};
    \silentremark{\sandor{060914}: I added ``$\exists\, C\subseteq X$'' here, because
      it seemed missing.\\ I also added here and elsewhere ``$c_1(\ )$'' to change a
      line bundle into a divisor...}
  \item  $\displaystyle{\mathbb{T}(L):=\bigcap_{\eps>0} T_\eps^A}$ to be the \emph{stable non-AEA locus of $L$}, and
  \item  $ ~ \mathbb{T}^0(L):=\set{x \with x \in C \textrm{ such that  }  L\cdot C < 0}$ the \emph{negative locus of $L$}.
  \end{enumerate}
\end{de}

\smallskip

\silentremark{\vskip-48pt\sandor{060914}: When labeling an ``item'', please use
  ``ilabel'' and refer to it with ``iref''. That will add the theorem number as well.
  (See sk-preamble.tex for definition).}

\begin{propdef}
\label{defline}
For  a line bundle $L$ on the  variety $X$ we have the following.
\begin{enumerate}

\item
\ilabel{nonnef-restricted}
$ \mathbb{T}^0(L) \subseteq \mathbb{B}_- (L)$, the inequality can be strict (cf.
\cite[Remark 6.3]{BDPP}).

\item
\ilabel{nonample-augmented}
 $ \mathbb{T}(L) \subseteq \mathbb{B}_+ (L)$.
%, the inequality can be strict  (cf. \sideremark{This needs a reference}).

\item $L$ is \emph{ample} iff $\mathbb{B}_+ (L) = \emptyset$.

\item $L$ is \emph{semiample} iff $\mathbb{B}(L) = \emptyset$.

\item $L$ is \emph{nef} iff $\mathbb{B}_- (L) = \emptyset$ iff $ \mathbb{T}^0(L) =
  \emptyset$.

\item
\ilabel{big}
$L$ is \emph{big} iff $\mathbb{B}_+ (L) \neq X$.

\item
\ilabel{psef}
$L$ is \emph{psef} (pseudo-effective) iff $\mathbb{B}_- (L) \neq X$.

\item
\ilabel{anef}
$L$ is \emph{almost nef} iff  $\mathbb{T}^0(L)$  is  contained
in a countable union of proper closed subsets of $X$.
%
%\item
%\ilabel{anef-psef}
%$L$ is almost nef if and only if it is pseudo-effective.
%
 \item
   \ilabel{aea}
   $L$ is \emph{AEA} (almost everywhere ample) if $ \mathbb{T} (L) \neq X$.

\item
\ilabel{wpos} $L$ is \emph{weakly positive} if $\overline{\mathbb{B}_- (L)} \neq X$.

\end{enumerate}

\end{propdef}

\begin{proof}

  Points (\ref{defline}.\ref{nonnef-restricted}-\ref{big}) are well-known statements.
  The claims \iref{anef}, \iref{aea} and \iref{wpos} are the definitions of respective notions
  according to \cite{BDPP}, \cite{Miy83} and \cite{vie}, respectively.  The only
  statement in need of a proof is \iref{psef}.

  Note that a line bundle is pseudo-effective precisely if  its numerical equivalence class lies in the closure of the
  effective cone in the real N\'eron-Severi group. Hence the line bundle  $L$ is psef if and only if  $\forall m>0
  ~ L+ (1/m) A$ is effective, or, equivalently, if  $\forall m>0 ~\mathbb{B} (mL + A) \neq X$. Therefore
  $\mathbb{B}_-(L) \neq X$ as it is contained in a countable union of proper closed
  subsets of $X$.  Conversely,  if $\mathbb{B}_-(L) \neq X$, then the class of $L$ is a
  limit of effective classes.
\end{proof}

\begin{remark}

  Positivity properties related to asymptotic base loci are best  summarized in the form of a table.

  \begin{center}

    \begin{tabular}{|c|c|c|c|c|}

      \hline

      & $\mathbb{B}_- (L)$ & ${\overline{\mathbb{B}_- (L)}  }$  & $ \mathbb{B} (L) $
      & $\mathbb{B}_+ (L)$ \vphantom{$\overline{B}^{\text{\normalsize B}}$}\\
      \hline
      $=\emptyset $& nef & nef &  semiample & ample \\
      \hline
      $\neq X $ & pseudo-effective & weakly positive & effective &  big \\

      \hline

    \end{tabular}

  \end{center}

\medskip

\end{remark}

\begin{remark}

  See Section~\ref{badaea} and in particular Remark \ref{remark-augm} for more
  details about non-AEA loci and their relationship with the augmented base loci.

  Furthermore, as we mentioned earlier,  Lesieutre \cite{lesieutre} proved
  that there exist line bundles which are pseudo-effective but not weakly positive. In
  particular, $\mathbb{B}_- (L)$ is not necessarily closed.
\end{remark}

\begin{proposition}

\label{anef-psef}
A line bundle $L$ is almost nef if and only if it is pseudo-effective.

\end{proposition}

\begin{proof}
  One implication is obvious by \iref{nonnef-restricted}.  The other implication
  follows from \cite{BDPP}, as if $L$ is not pseudo-effective then there exists a
  reduced irreducible curve $C \subseteq X$, such that $c_1(L)\cdot C <0$ and $C$
  moves in a family covering all $X$, so $ \mathbb{T}^0(L)$ cannot be contained in a
  countable union of proper (Zariski) closed subsets.
\end{proof}

The following theorem will be proved in Section~\ref{sectionaea}:

\begin{theorem}

A line bundle $L$ is big if and only if it is AEA.

\end{theorem}

\silentremark{\sandor{060914}: I changed ``AEN'' to ``AEA'' assuming that it was a
  typo.}

A recent result of Lehmann \cite{brian} gives a characterization of the relationship
between the non-AEA locus and the diminished base locus. We will use the following
when describing all the relationships.

\begin{de} Let $X$ be a smooth projective variety over $\mathbb{C}$ and let $D$ be a
  pseudo-effective $\R$-divisor on $X$. Suppose that $\phi : Y \to X$ is a proper
  birational map from a smooth variety $Y$ . The movable transform of $L$ on $Y$ is
  defined to be
  $$
  \phi_{\mov}^{-1}(D) := \phi^*D - \sum_{E ~ \phi-exc} \sigma_E(\phi^*D)\cdot E.
  $$
  For a pseudo-effective line bundle $L$  define $\phi_{\mov}^{-1}(L)$ as the line
  bundle associated to $\phi_{\mov}^{-1}(c_1(L))$ on $Y$.
\end{de}

Note that the movable transform is not linear and is only defined for
pseudo-effective divisors.

\begin{remark}
  In the above, $\sigma_E$ is the asymptotic multiplicity function introduced by
  Nakayama \cite[Section III.1]{Nak04}.
  If $X$ is a smooth projective variety, $L$ a pseudo-effective $\R$-divisor, $E$ a
  prime divisor on $X$, then
  \[
  \sigma_E(L) := \lim_{\epsilon\to 0^+} \inf \{ \mult_E{L'}\mid L' \geq 0 \text{ and
  } L'\sim_{\R} L+\epsilon A\}\ ,
  \]
  where $A$ is an arbitrary but fixed ample divisor.
\end{remark}

\begin{remark}
  Following \cite[Definition 1.2]{brian}, we call an irreducible curve $C$ on $X$ to
  be a $\mov^1$-curve, if it deforms to cover a codimension one subset of $X$.
\end{remark}

\begin{thm}[\cite{brian}] Let $X$ be a smooth projective variety over $\C$ and $D$ a
  pseudo-effective $\R$-divisor. $D$ is not movable if and only if there is a
  $\mov^1$- curve $C$ on $X$ and a proper birational morphism $\phi: Y \to X$ from a
  smooth variety $Y$ such that $$\phi_{\mov}^{-1}(D) \cdot \bar{C} < 0,$$ where
  $\bar{C}$ is the strict transform of a generic deformation of $C$.
\end{thm}

The following reformulation is easy to see.

\begin{thm}
  Let $X$ be a smooth projective variety and $D$ a pseudo-effective $\R$-divisor.
  Suppose that $V$ is an irreducible subvariety of $X$ contained in $\mathbb{B}_-(L)$
  and let $\psi : X' \to X$ be a smooth
  birational model resolving the ideal sheaf of $V$.  Then there is a birational
  morphism $\phi: Y \to X'$ from a smooth variety $Y$ and an irreducible curve
  $\bar{C}$ on $Y$ such that $\phi^{-1}_{\mov}(\psi^*D)\cdot \bar{C} < 0$ and $\psi
  \circ \varphi (\bar{C})$ deforms to cover $V$ .

\end{thm}

\begin{remark} Let $L$ be any line bundle on $X$ smooth projective, then
$$\mathbb{B}_-(L)= \bigcup_{f: Y \to X} f(\mathbb{T}^0(f_{\mov}^{-1}(L))).$$

\end{remark}

There are several examples for which the loci $\mathbb{T}^0(L)$ and $\mathbb{B}_-(L)$
do not coincide, and in some cases the difference is divisorial.

\begin{qsn}
Is it true that $\mathbb{B}_-(L)$ is contained in a proper closed subset of $X$
(i.e., $L$ is weakly positive) if and only if the same holds for $\mathbb{T}^0(L)$?
\end{qsn}

\begin{remark}
In \cite[Question 7.5]{BDPP} the authors ask if for a vector bundle $E$, $\mathbb{B}_-(\sO_{\mathbb{P}(E)}(1))$
doesn't dominate $X$ if and only if neither does $\mathbb{T}^0(\sO_{\mathbb{P}(E)}(1))$.
\end{remark}

\section{Restricted base loci for vector bundles}
\label{restricted}

Here we explore the connections between  the positivity properties of a vector bundle
and  the associated asymptotic base loci.  We will start recalling some classical
definitions for vector bundles. Note that these definitions do sometimes appear
slightly differently in the literature, but we will try to follow and indicate specific
selected references each time.

\begin{definition}

  Let $E$ be a vector bundle on a smooth projective variety $X$, $\pi \colon
  \mathbb{P}(E) \to X$ the projective bundle of one dimensional quotients of $E$, and
  $\sO_{\mathbb{P}(E)}(1)$ the universal quotient of $\pi^* E$ on $\mathbb{P}(E)$.
  We say that $E$ is

  \begin{enumerate}

  \item \emph{nef} if $\sO_{\mathbb{P}(E)}(1))$ is a nef line bundle, {i.e.,} if
    $\mathbb{B}_- (\sO_{\mathbb{P}(E)}(1)) = \emptyset$;

  \item \emph{almost nef} if $\pi (\mathbb{T}^0(\sO_{\mathbb{P}(E)}(1)))$ is
    contained in a countable union of proper closed subsets of $X$ (cf. \cite{BDPP});

  \item \emph{pseudo-effective} if $\mathbb{B}_- (E) \neq X$ (cf. \cite{BDPP});

  \item \emph{weakly positive} if $\overline{\mathbb{B}_- (E)} \neq X$ (cf.
    \cite{vie}); \ilabel{wp}

  \end{enumerate}

\end{definition}

\begin{proposition}

A vector bundle $E$ is nef if and only if $\mathbb{B}_- (E)= \emptyset$.

\end{proposition}

\begin{proof}

This follows directly from Proposition~\ref{bmin}.
\end{proof}

\begin{proposition}

  A vector bundle $E$ is pseudo-effective if and only if $\sO_{\mathbb{P}(E)}(1)$ is
  pseudo-effective and $\pi (\mathbb{B}_- (\sO_{\mathbb{P}(E)}(1))) \neq X$.

\end{proposition}

\begin{proof}
  This follows immediately from Proposition~\ref{bmin}.
\end{proof}

\begin{remark}

  The proposition above is the same as \cite[Proposition 7.2]{BDPP}. Observe that the
  locus $L_{\nonnef}$ in \cite{BDPP} is what we call $\mathbb{B}_- (L)$ here (we are
  in the smooth projective case).

\end{remark}

The following proposition is immediate.

\begin{proposition}

  A vector bundle $E$ is almost nef if and only if there exists a countable union $T
  = \bigcup_{i \in \mathbb{N}} T_i$ of proper closed subsets of $X$, such that for
  any curve $C \subseteq X$ not contained in $T$ the restriction $E|_{C}$ is a nef
  vector bundle.

\end{proposition}

\begin{remark}

  It follows from the definitions and propositions above that
  \[
  E ~\textrm{weakly positive} \Rightarrow E ~\textrm{pseudo-effective} \Rightarrow E
  ~\textrm{almost nef} ~.
  \]
  We have seen that the first implication is not an equivalence, while it is an open
  question whether, for vector bundles, being almost nef is equivalent to being
  pseudo-effective, as in the line bundle case cf.\ \cite[Question~7.5]{BDPP}.

\end{remark}

\begin{question}
Does $E$ being {almost nef} imply that $E$ is {pseudo-effective}?
\end{question}

If $E$ is almost nef, then the line bundle $\sO_{\mathbb{P}(E)}(1)$ is almost nef,
hence pseudo-effective.  In order to have that $E$ is pseudo-effective, one needs to
show that $\mathbb{B}_- (\sO_{\mathbb{P}(E)}(1))$ does not dominate $X$.

%% Does \cite{brian} imply that $E$ is also pseudo-effective?  Or, on the contrary,
%% can we provide examples were $E$ is almost nef but for example L-big and not
%% pseudo-effective?

\section{Augmented base loci for vector bundles}
\label{augmented}

\begin{definition}

  Let $E$ be a vector bundle on a smooth projective variety $X$, $\pi \colon
  \mathbb{P}(E) \to X$ the projective bundle of one dimensional quotients of $E$, and
  $\sO_{\mathbb{P}(E)}(1)$ the universal quotient of $\pi^* E$ on $\mathbb{P}(E)$.

\begin{enumerate} We say that

\item $E$ is \emph{ample} if $\sO_{\mathbb{P}(E)}(1)$ is ample on $\mathbb{P}(E)$;

\item $E$ is \emph{L-big} if $\sO_{\mathbb{P}(E)}(1)$ is big on $\mathbb{P}(E)$; and

\item $E$ is \emph{V-big} (or \emph{Viehweg-big}) if there exists an ample line
  bundle $A$ and a positive integer $c>0$ such that $\Sym^c E \otimes A^{-1}$ is
  weakly positive, {i.e.,} such that $\overline{\mathbb{B}_- (\Sym^q E \otimes
    A^{-1})} \subsetneq X$ (cf.\ \iref{wp}).

\end{enumerate}

\end{definition}

\begin{proposition}
  A vector bundle $E$ is ample if and only if $\mathbb{B}_+ (E)= \emptyset$.
\end{proposition}

\begin{proof}
Follows directly from Proposition~\ref{bplus}.
\end{proof}

\begin{remark}

  It is well-known that a line bundle $L$ is big if and only if $\mathbb{B}_+ (L)
  \neq X$, equivalently,  if there exist $A$ ample and a positive integer $c>0$ such that
  $L^{\otimes c} \otimes A^{-1}$ is pseudo-effective.  Next we will prove that the
  same equivalences hold for V-bigness for vector bundles of arbitrary rank.

\end{remark}

\begin{theorem}
  \label{thm:V-big}

  Let $E$ be a vector bundle on a smooth projective variety $X$.  Then the following
  properties are equivalent:

\begin{enumerate}

\item $E$ is V-big.
\ilabel{one}

\item There exist an ample line bundle $A$ and a positive integer $c>0$ such that
$\Sym^c E \otimes A^{-1}$ is pseudo-effective.
\ilabel{two}

\item $\mathbb{B}_+ (E) \neq X$.
\ilabel{three}

\end{enumerate}

\end{theorem}

\begin{proof}

  The implication $\iref{one} \Rightarrow \iref{two}$ is clear; let us consider
  $\iref{two} \Rightarrow \iref{three}$.  Suppose there exist $A$ ample and $c>0$
  such that $\Sym^c E \otimes A^{-1}$ is pseudo-effective, {i.e.,}
 \[
 \mathbb{B}_-(\Sym^c E \otimes A^{-1}) \neq X \ .
 \]
 Then
\[
\bigcup_{q>0} \mathbb{B} (E -(1/c)A + (1/q)A) =
\bigcup_{q>0} \mathbb{B} (E - \frac{q-c}{qc}A) \neq X ~.
\]
Now for $q\gg0$ and $\frac{q-c}{qc} > \frac{1}{2c}$, one has
\[
\mathbb{B}_+ (E) \,\subseteq\, \mathbb{B} (E - \frac{1}{2c}A)
\,\subseteq\, \mathbb{B} (E - \frac{q-c}{qc}A) \,\subsetneq\, X  \ ,
\]
hence the validity of the desired implication.

Next we verify that $\iref{three} \Rightarrow \iref{one}$.
If $E$ satisfies $\mathbb{B}_+ (E) \neq X$,
then there exists $q>0$ such that $\mathbb{B} (E - \frac{1}{q}A) \subsetneq X$ is a closed proper subset.
Consequently,

\[
\overline{\mathbb{B}_- (\Sym^q E \otimes A^{-1})} \,\subseteq\, \mathbb{B} (\Sym^q E
\otimes A^{-1}) \,\subsetneq\, X \qedhere
\]
\end{proof}

\begin{corollary}
  If $E$ is V-big, then it is  L-big as well.
\end{corollary}

\begin{proof}
  Theorem~\ref{thm:V-big} yields $\mathbb{B}_+ (E) \neq X$, therefore
  $\mathbb{B}_+ (\sO_{\mathbb P(E)}(1)) \neq \mathbb P(E)$ via Proposition~\ref{bplus}. As a consequence
  $\sO_{\mathbb P(E)}(1)$ is big on $\mathbb P(E)$, equivalently,  $E$ is L-big.
\end{proof}

\begin{remark}\label{L-big-vs-V-big}
  L-big vector bundles are not necessarily V-big, as the example of
  $\sO_{\P^1}\oplus\sO_{\P^1}(1)$ on $\P^1$ shows (see \cite[p.24]{KellyThesis}).

  The key difference between V-big and L-big vector bundles is that being L-big means
  that $\sO_{\P(\she)}(1)$ is ample with respect to an open set $V\subseteq \P$
  whereas $\she$ is V-big if we can take $V$ to be of the form $V=\pi^{-1}(U)$,
  $U\subseteq X$ open.
\end{remark}

\begin{remark}

  A vector bundle $E$ on a variety $X$ satisfying $\mathbb{B}_+ (E) \neq X$ is also
  called \emph{ample with respect to an open subset} (cf. \cite[Chapter
  3]{KellyThesis}).

\sideremark{\sandor{060914}: this should be properly defined. I'll do this later.}

\end{remark}

\begin{remark}
\label{diverio}

In the case where $E = \Omega_X$ is the cotangent sheaf of a variety $X$, the
definitions and Proposition~\ref{bplus} imply the following inclusion $\mathbb{B}_+
(\Omega_X) \supseteq DS(X, T_X)$, where $DS(X, T_X)$ is the Demailly-Semple locus.

The work of Diverio and Rousseau \cite{divrous} therefore provides examples of
complex projective varieties of general type $X$ where $\Omega_X$ is a semistable
L-big vector bundle with a big determinant, which is nevertheless not V-big.

\end{remark}

\section{Almost Everywhere Ampleness}
\label{sectionaea}

The notion of almost everywhere ampleness was first defined by Miyaoka in the context of his work on vector bundles on surfaces; the definition
goes through in all dimensions verbatim.

\begin{definition}[\cite{Miy83}]
  Let $X$ be a smooth projective variety,  $E$ a rank $r$ vector bundle on $X$.
  Consider the projectivized bundle $\P=\P(E)$ with projection morphism $\pi:\P\to X$ and
  tautological bundle $\sO_\P(1)$. We say that $E$ is \emph{almost everywhere ample}
  (AEA for short), if there exists an ample line bundle $A$ on $X$, a Zariski closed
  subset $T\subset \P$, whose projection $\pi(T)$ onto $X$ satisfies $\pi(T)\neq X$,
  and a positive number $\eps>0$ such that
   $$
   c_1(\sO_\P(1))\cdot C\ge \eps\cdot \pi^*(c_1(A))\cdot C
   $$
   for all curves $C\subset \P$ that are not contained in $T$.
\end{definition}

\begin{comment}
  It is possible that the scale-invariant version will turn out to be a more useful
  definition: $E$ is AEA if the condition in the definition is fulfilled for
  $\Sym^mE$ for some $m>0$.  This changes nothing about our results for line bundles,
  but might give a better Kodaira lemma.
\end{comment}

For line bundles, this notion coincides with bigness:

\begin{proposition}\label{proposition-AEA-big}
   For a line bundle $L$ on a smooth projective variety $X$, the
   following are equivalent:
   \begin{enumerate}
   \item $L$ is AEA, i.e, there is an ample line bundle $A$ on $X$, a number
     $\eps>0$, and a proper Zariski closed subset $T\subset \P$ such that
     $$
     c_1(L)\cdot C\ge \eps\cdot c_1(A)\cdot C
     $$
     for all curves $C\subset \P$ not contained in $T$.
     \ilabel{onec}
   \item For every ample line bundle $A$ on $X$, there is an $\eps>0$ and a proper
     Zariski closed subset $T\subset X$ such that
     $$
     c_1(L)\cdot C\ge \eps\cdot c_1(A)\cdot C
     $$
     for all curves $C\subset X$ not contained in $T$.
     \ilabel{twoc}
   \item
     $L$ is big.
     \ilabel{threec}
   \end{enumerate}
\end{proposition}

\begin{proof}
  Assume \iref{threec}, and let $A$ be any ample line bundle.  Then, by Kodaira's lemma,
  there is a positive integer $m$ such that we can write
   $$
   mc_1(L) = c_1(A)+F ,
   $$
   where $F$ is an effective divisor.  Taking $T$ to be the support of $F$, it
   follows for every curve $C\subset X$ not contained in $T$ that
   $$
   mc_1(L)\cdot C=c_1(A)\cdot C+F\cdot C\ge c_1(A)\cdot C
   $$
   and this implies \ref{twoc} with $\eps:=1/m$.

   Obviously \ref{twoc} implies \iref{onec}, so let us assume condition \iref{onec}
   and show that it implies \iref{threec}.  A curve $C\subset X$ such that $c_1(L)\cdot
   C<\eps\cdot c_1(A)\cdot C$ cannot be a movable curve (in the sense of
   \cite[Sect.~11.4.C]{PAG}), since these cover all of $X$ (by \cite[Lemma
   11.4.18]{PAG}), whereas $T\ne X$.  So $L$ must have positive intersection with all
   movable curves.  This implies that $L$ lies in the dual of the cone of movable
   curves $\Mov(X)$, which by the theorem of Boucksom-Demailly-Paun-Peternell
   \cite{BDPP} is the pseudo-effective cone $\overline{\Eff}(X)$.  In order to
   conclude that $L$ is big -- and thus to complete the proof -- it is therefore
   enough to show that $L$ lies in the interior of that cone.

   The assumption that $L$ be AEA says that
   $$
   (c_1(L)-\eps c_1(A))C\ge 0 \qquad\mbox{for all } C\not\subset T \ .
   $$
   Therefore, writing $c_1(L)-\eps c_1(A)=(c_1(L)-\frac\eps2 c_1(A))-\frac\eps2
   c_1(A)$, we see that $c_1(L)-\frac\eps2 c_1(A)$ is AEA as well. Moreover, every
   class in the open set
   $$
   c_1(L)-\frac\eps2 c_1(A)+\Amp(X)
   $$
   is AEA, and $c_1(L)$ lies in this open set.
\end{proof}

\begin{proposition}\label{proposition-AEA}
  Let $E$ be a vector bundle on a smooth projective variety $X$, let $\P=\P(E)$. If
  $E$ is AEA on $X$, then so is $\sO_\P(1)$ on $\P$.
\end{proposition}
\begin{proof}
  For $E$ to be AEA means that for every ample line bundle $A$ on $X$, there exists a
  Zariski-closed subset $T\subseteq X$, and $\eps >0$ such that
  \[
  c_1(\sO_\P(1))\cdot C \geq \eps (\pi^*c_1(A)\cdot C)
  \]
  for all irreducible curves not contained in $T$.

  Since $\sO_\P(1)$ is $\pi$-ample, the line bundle $\pi^*A\otimes \sO_\P(m)$ is
  ample for all $m\geq m_0\gg 0$ by \cite[Proposition 1.7.10]{PAG}.  According to
  Proposition~\ref{proposition-AEA-big}, $\sO_\P(1)$ is AEA if and only if it is big,
  therefore it suffices to prove the AEA property for $\sO_\P(k)$ for some large $k$.
  This means in particular, that we are allowed to work with $\Q$-divisors as well.

  Let $m_0$ be, as above, a positive integer such that $\pi^*A\otimes\sO_\P(m_0)$ is
  ample.  We will prove that $\sO_\P(1)$ is AEA on $\P$ with closed subset
  $T\nsubseteq \P$, and a suitable $\eps'>0$.  We need that
  \[
  c_1(\sO_\P(1))\cdot C \geq \eps'(c_1(\pi^*A\otimes\sO_\P(m_0))\cdot C) \ ,
  \]
  or equivalently,
  \[
  c_1(\sO_\P(1))\cdot C \geq \frac{\eps'}{1-\eps'm_0}(\pi^*c_1(A)\cdot C)
  \]
  for all curves not contained in $T$. By our assumption on $E$, this holds whenever
  \[
  \eps' < \frac{\eps}{1+\eps m_0}\ . \qedhere
  \]
\end{proof}

\begin{corollary}
  If $E$ is AEA, then it is L-big.
\end{corollary}
\begin{proof}
  Immediate from Proposition~\ref{proposition-AEA-big} and
  Proposition~\ref{proposition-AEA}.
\end{proof}

%
%
%\begin{remark}
%  My impression is that the converse does not hold, it rather seems the AEAness is
%  equivalent to ampleness with respect to an open subset. The counterexample
%  $\sO\oplus\sO(1)$ on $\P^1$ showing that not all big vector bundles are ample with
%  respect to an open set should work here as well.
%\end{remark}

%*****************************************************************************

\section{The bad AEA locus in the line bundle case}
\label{badaea}

Consider a line bundle $L$, and fix an ample line bundle $A$ and a number $\eps>0$.
We defined the \textit{non-AEA locus} of $L$ with respect to $A$ and $\eps$ as the
subvariety
$$
T_\eps^A=\mbox{closure}(\bigcup\set{C\with C\mbox{ curve on $X$ with $c_1(L)\cdot
    C<\eps\cdot c_1(A)\cdot C$}})\ .
$$
The AEA assumption on $L$ simply  means that there exists an $\eps>0$ such that
$T^A_\eps\ne X$.  For $\eps<\delta$ we have $T_\eps\subset T_\delta$, so that we can
 express the AEA condition equivalently as saying that the intersection
${\mathbb{T}}(L):=\bigcap_{\eps>0} T_\eps^A$ is not all of $X$.

%
%\begin{remark}\label{remark-char-T}
%   We have
%   $$
%      T(L)=\bigcup\set{C\with c_1(L)\cdot C\le 0} \ .
%   $$
%   In fact, if $C$ is a curve with $c_1(L)\cdot C\le 0$, then obviously
%   it satisfies $c_1(L)\cdot C\le \eps c_1(A)\cdot C$ for every $\eps>0$, so
%   that is lies in every $T_\eps$, and hence in $T(L)$.
%   If, conversely, $c_1(L)\cdot C\le \eps c_1(A)\cdot C$
%   holds for every $\eps>0$, then
%   we must have $c_1(L)\cdot C\le 0$, so $C$ lies in $T(L)$.
%\end{remark}

\begin{remark}\label{remark-augm}
   It is immediate that
   $$
       \mathbb{T}(L)\subseteq \mathbb{B}_+(L) \ .
   $$
   In fact, by the noetherian property there are positive real numbers $\eps_0$ and
   $\delta_0$ such that

   \[
   \mathbb{T}(L) = T_{\eps}^A ~~\forall \eps \leqslant \eps_0 \textrm{ and }
   \mathbb{B}_+ (L) = \mathbb{B}(L - \delta A) ~~\forall \delta \leqslant \delta_0.
   \]

   Now choose $\eps < \min (\eps_0 , \delta_0)$, then
   \[
   \mathbb{T}(L) = T_{\eps} = \overline{\set{ x \with x \in C \mbox{ curve on $X$
         s.t. } c_1(L)\cdot C<\eps\cdot c_1(A)\cdot C}}.
   \]
   If $C$ is a curve such that $(L - \eps A)C <0$ then $C \subseteq \mathbb{B}(L -
   \eps A) = \mathbb{B}_+ (L)$ which is a closed set.

   %
   % if $C$ is a curve with $c_1(L)\cdot C\le 0$, then for every $\eps > 0$ we have
   % $$
   % (L-\eps A)C<0
   % $$
   % which implies that $C$ is contained in the stable base locus of $L-\eps A$. But
   % for small $\eps$ the latter equals $B_+(L)$.
\end{remark}

\begin{remark}

  A line bundle $L$ is ample if and only if $ \mathbb{T}(L) = \emptyset$. In general
  the inclusion
  \silentremark{\sandor{060914}: should we mention Seshadri constants here?}
  \[
  \overline{\set{ x \with x \in C \mbox{ curve on $X$ s.t. } c_1(L)\cdot C \leqslant
      0}} \subseteq \mathbb{T}(L)
  \]
  is strict, as shown by a strictly nef non ample line bundle $L$, where the first
  set is empty but the second one is not.  Examples of line bundles that are strictly
  nef (and even big) and non ample have been first given by Mumford (cf.
  \cite{har70}), and a complete description can be found in \cite{urbin}.

\end{remark}

\begin{remark} A few words on the relationship between $\mathbb{T}(L)$ and $\mathbb{B}_\pm(L)$.
  We'll show here that $\mathbb{T}(L)\neq B_+(L)$ in general. A bit more precisely,
  we will try to understand the relationship of $\mathbb{T}(L)$ to the augmented and
  restricted base loci of $L$ when $\dim X=2$. Recall that
\[
 \mathbb{B}_-(L) = \bigcup_{m=1}^{\infty} \mathbb{B} (L+\frac{1}{m}A)
\]
for any integral ample divisor $A$ on $X$.

Let $D$ be a big divisor on a smooth projective surface $X$ with Zariski
decomposition $D=P_D+N_D$. Then \cite[Examples 1.11 and 1.17]{ELMNP} tell us that
\begin{eqnarray*}
  \mathbb{B}_+(\sO_X(D)) & = & \Null(P_D) \equ \bigcup_{C}\st{C\subseteq X \text{ irred
    }\,|\, P_D\cdot C=0}\ , \\
  \mathbb{B}_-(\sO_X(D)) & = & \Supp N_D\ .
\end{eqnarray*}
\end{remark}

\begin{example} Here we present an  example where $T(L)\neq \mathbb{B}_{\pm}(L)$.
  Let $X$ be a surface that carries a big divisor $D$ and an irreducible curve
  $C\subseteq X$ satisfying $C\subseteq \Supp N_D$ and $D\cdot C>0$. Then
  $C\not\subseteq T(L)$, but $C\subseteq \mathbb{B}_-(\sO_X(D))\subseteq
  \mathbb{B}_+(\sO_X(D))$.  In this case we have
   \[
   T(\sO_X(D)) \neq \mathbb{B}_-(\sO_X(D)),\mathbb{B}_+(\sO_X(D)).
   \]
   Surfaces carrying such $D$ and $C$ exist by \cite{BauFun10}: Consider a K3 surface
   $X$, on which the Zariski chamber decomposition does not coincide with the Weyl
   chamber decomposition. The latter is by \cite[Theorem~1]{BauFun10} the case if and
   only if there are $(-2)$-curves on $X$ having intersection number 1. For a
   concrete example one can, as done in \cite[Section~3]{BauFun10}, take a smooth
   quartic surface $X\subset\mathbb P^3$ that has a hyperplane section of the form
   $H=L_1+L_2+Q$, where $L_1$ and $L_2$ are lines and $Q$ is a smooth conic. Then the
   divisors of the form
   $$
      D=H+a_1L_1+a_2L_2
   $$
   with $a_1\ge 1$ and $a_2\ge 1$ have $L_1$ and $L_2$ in the support of the negative
   part of their Zariski decomposition, but one can find $a_1,a_2$ such that $D\cdot
   L_1>0$ and $D\cdot L_2<0$ (for instance $a_1=2,a_2=4$). (In the notation of
   \cite{BauFun10}, $D$ lies in the Zariski chamber $Z_{\{L_1,L_2\}}$, but in the
   Weyl chamber $W_{\{L_2\}}$.)
\end{example}

\begin{example}
  In general $T(\sO_X(D))$ is not contained in $B_-(\sO_X(D))$ either, where $L=\sO_X(D)$ for a suitable integral Cartier  divisor $D$.
  To see this,
  take a surface where all negative curves have self-intersection $-1$. Then the
  intersection form of the negative part of the Zariski decomposition of any big
  divisor is $-\Id$, in other words, no two curves in it can intersect. Consequently,
  \[
  T(\sO_X(D)) \equ B_+(\sO_X(D))\ .
  \]
  This can be seen as follows: let $C\subseteq B_+(\sO_X(D))=\Null (P_D)$ be an
  irreducible curve. Since $P_D$ is big and nef, the intersection form on $P_D^\perp$
  is negative definite, which under the given circumstances means that $(C\cdot
  C')=0$ for every irreducible curve $C\neq C'$ coming up in $N_D$. Therefore
  \[
  D\cdot C \equ P_D\cdot C + N_D\cdot C \equ 0 + (\leq 0) \equ (\leq 0)\ .
  \]
  Consequently, $C\subseteq T(\sO_X(D))$.

  Take a non-stable (in the sense if \cite[Definition 1.22]{ELMNP} big divisor $D$ on
  $X$, then $B_-(\sO_X(D))\subsetneq B_+(\sO_X(D)) = T(\sO_X(D))$.
\end{example}

\medskip
The following lemma is well-known to experts working in the area, but for lack of an
adequate reference we include it here.

\begin{lemma}
  Let $D$ be a big divisor on a smooth projective surface, $C\subseteq X$ irreducible
  curve, $D\cdot C=0$. Then $(C^2)<0$.
\end{lemma}
\begin{proof}
  Let $D=P_D+N_D$ denote the Zariski decomposition of $D$. If $C\subseteq \Supp N_D$,
  then it must have negative self-intersection, since the intersection form on $N_D$
  is negative definite. Assume $C$ is not in $N_D$.  Then
  \[
  0 \equ D\cdot C \equ P_D\cdot C +N_D\cdot C \equ (\geq 0) + (\geq 0) \ ,
  \]
  since $P_D$ is nef, $C$ is effective, $N_D$ is effective with no common components
  with $C$. This can only happen if
  \[
  P_D\cdot C \equ N_D\cdot C \equ 0\ .
  \]
  Therefore, $C$ is orthogonal to the big and nef divisor $P_D$, hence we must
  have $(C^2)<0$.
\end{proof}

\section{V-big vs. AEA}

Let $X$ be a smooth projective variety and $\she$ a vector bundle on $X$.  There
exist two non-equivalent definitions for bigness in the literature: V-big and L-big
vector bundles.  It is known that V-bigness implies L-bigness and that the converse
does not hold if $\rk\she\geq 2$ (cf.\ Remark~\ref{diverio}).

Throughout this section we will point out some differences (for example a different
Kodaira's lemma) between L-big and V-big vector bundles, and compare V-bigness and
almost everywhere amplenessq.  In particular we will show that these positivity
properties coincide for nef vector bundles.  V-big vector bundles are also called
\emph{ample with respect to an open set} \cite[Chapter 3]{KellyThesis}.
% , so we will refer to them in such way when stating properties proven in
% \cite{KellyThesis}.

We have seen that if $E$ is V-big, then $E$ is also AEA cf.\
Remark~\ref{remark-augm}.

\begin{question}
  Does $E$ being AEA imply that $E$ is V-big?

\noindent
We will see that this is the case if $E$ is nef.

\end{question}

\begin{remark}
  As pointed out in \cite[Lemma 3.44]{KellyThesis}, a vector bundle on a projective
  curve is ample with respect to an open set exactly if it is ample.
\end{remark}

Next we will show that a strong form of Kodaira's lemma is valid for vector bundles
that are ample with respect to an open set.

\begin{lemma}[(Kodaira's lemma for vector bundles)]
  Let $X$ be a smooth projective variety, $\she$ a vector bundle, and $\sha$ an ample
  line bundle on $X$. Then the following are equivalent.
  \begin{enumerate}
  \item $\she$ is ample with respect to an open subset.
    \ilabel{oneb}
  \item $\Sym^m\she$ contains an ample vector bundle of the same rank for some $m>0$.
    \ilabel{twob}
  \item There exists $m>0$ and an injective morphism
    \[
    \bigoplus^{\rk\Sym^m\she}\sha \hookrightarrow \Sym^m\she\ ,
    \]
    which is an isomorphism over an open subset.
    \ilabel{threeb}
\end{enumerate}
\end{lemma}
\begin{proof}
  The equivalence of \iref{oneb} and \iref{threeb} is the content of \cite[Lemma
  3.42]{KellyThesis}; \iref{threeb} obviously implies \iref{twob}, and
  \iref{twob} implies \iref{oneb} holds if being ample with respect to an open set is
  scale-invariant.
  \sideremark{\sandor{060914}: this needs to be rewritten after including the defn of
    ample wrt an open set}
\end{proof}

There is a useful characterization of ampleness with respect to an open subset in terms of $\sO_{\P(\she)}(1)$.

\begin{lemma}(cf.\ Proposition~\ref{bplus})\label{lemma-equiv} With notation as
  above, $\she$ is ample with respect to the dense open set $U\subseteq X$ precisely
  if $\sO_{\P(\she)}(1)$ is ample with respect to $\pi^{-1}(U)\subseteq \P(\she)$.
\end{lemma}

\begin{comment}
 With notation as above, recall that  a vector bundle $\she$ is called \emph{L-big},
 if the line bundle $\sO_{\P(\she)}(1)$ on $\P(\she)$ is big.
\end{comment}
\medskip

Here we have the following weaker version of the Kodaira lemma.

\begin{lemma}
  % (\cite[Example 6.1.23]{PAG})
  \label{lemma-PAG}
  Let $\she$ be a vector bundle using the notation above.
\begin{enumerate}

\item Assume that $\HH{0}{X}{\Sym^m\she}\neq 0$ for some $m>0$.  Then for any ample
  line bundle $\sha$ on $X$ and any $k>0$, ${\Sym^k\she}\otimes\sha$ is L-big.
  \ilabel{oned}

\item Assume that for some $m>0$ and some $x \in X$ the vector bundle $\Sym^m\she$ is
  generated at $x$ by its global sections $\HH{0}{X}{\Sym^m\she}$.  Then for any
  ample line bundle $\sha$ on $X$ and any $k>0$, ${\Sym^k\she}\otimes\sha$ is V-big.
  \ilabel{twod}

\item Conversely, assume that $\she$ is L-big. Then for any line bundle $\shl$ on
  $X$,
  \[
  \HH{0}{X}{\Sym^m\she\otimes \shl}\neq 0
  \]
  for all $m\gg 0$.
  \ilabel{threed}

\item Assume that $\she$ is V-big. Then for any line bundle $\shl$ on $X$, \(
  {\Sym^m\she\otimes \shl} \) is generically generated by its global sections for all
  $m\gg 0$.
  \ilabel{fourd}

\end{enumerate}

\end{lemma}

\begin{proof}
  To prove \iref{oned}, assume that $\HH{0}{X}{\Sym^m\she}\neq 0$ for some $m>0$.
  This means that $\HH{0}{\P(\she)}{\sO_{\P(\she)}(m)}\simeq \HH{0}{X}{\Sym^m\she}
  \neq 0$, hence $\sO_{\P(\she)}(1)$ is $\Q$-effective. Then
  \[
  \sO_{\P(\she\otimes\sha)}(1)\simeq \sO_{\P(\she)}(1)\otimes\pi^*\sha\ .
  \]
  By \cite[Proposition 1.7.10]{PAG}, the $\Q$-divisor $a c_1(\op) +
  \pi^*c_1(\shl)$ is ample for $0<a\ll 1$. This implies that
  \[
  c_1(\op)+\pi^*c_1(\shl) \equ (1-a) c_1(\op) + (a c_1(\op)+\pi^*c_1(\shl))
  \]
  can be written as the sum of an effective and an ample divisor, hence it is big.

  To prove \iref{twod}, assume that for some $m>0$ and some $x \in X$ the vector
  bundle $\Sym^m\she$ is generated at $x$ by its global sections.
  % $\HH{0}{X}{\Sym^m\she}$.
  Then $\mathbb{B}({\she}) \neq X$ and hence $\mathbb{B}_- ({\she}) \subset
  \mathbb{B}({\she}) \neq X$. Therefore ${\she}$ is weakly positive and for all $H$
  ample \[\mathbb{B}_+ ({\she} + (1/m) H) = \bigcap \mathbb{B}({\she} + (1/m)H -
  (1/n)H) \subseteq \mathbb{B}({\she}) \neq X.\]

  \iref{threed} is a reformulation of the Kodaira lemma on $\P(\she)$ (see
  \cite[Lemma 2.2.6]{PAG}).

  Finally, for \iref{fourd}, assume that $\she$ is V-big, $H$ an ample line bundle,
  and $\shl$ any line bundle on $X$.  Using notations of the first chapters we have
  that for $m$ sufficiently large $H - (1/m)L$ is ample, and as ${\she}$ is V-big we
  have
  \[
  \mathbb{B}(\Sym^{mk}{\she} +L) = \mathbb{B}(\Sym^k {\she} + (1/m) L) \subseteq
  \mathbb{B} (\Sym^k {\she} - (H - (1/m) L)) \neq X
  \]
  for $k$ big enough.
\end{proof}

\begin{remark}
  A line bundle is (L-)big if and only if it is ample with respect to an open set.
  We have seen in \ref{L-big-vs-V-big} that there exist simple vector bundles that
  are L-big, but not V-big, and hence not ample with respect to an open set.
\end{remark}

\begin{comment}
  It follows from Lemma~\ref{lemma-equiv} and the previous remark that a vector
  bundle, which is ample with respect to an open set is necessarily L-big.
\end{comment}

%\begin{framed}
%\end{framed}

\begin{proposition}
  Let $X$ be a smooth projective variety, $\she$ a vector bundle on $X$ containing an
  ample vector bundle $\sha$ of the same rank. Then
  \begin{enumerate}
  \item $\she$ is almost everywhere ample with respect to the closed subset $T= \Supp
    ({\she}/A)$.
  \item $\she$ is ample with respect to the dense open set $X\setminus T$.
  \end{enumerate}
\end{proposition}

\sideremark{\sandor{060914}: This proof is a mess. I will clean this up after adding
  the def of ample wrt an open set.}

\begin{proof}
%  {\color{gray}{
  Let $\shh$ be an ample line bundle on $X$, and $\shl=\sO_{\P(\she)}(1)$ and
  $C\subseteq \P(\she)$ an irreducible curve. If $C$ is contained in a fibre that
  maps to a point away from $T$, then
  \[
  \shl \cdot C \equ \pi^*\sha\cdot C\ ,
  \]
  which is good. If $C$ is contained in a fibre mapping to $T$, then we do not care
  about the intersections numbers at all.

  We may now assume that $C$ is not contained in a fibre of $\pi$. Let $B\deq
  \pi(C)\subseteq X$.

  By restricting everything to $B$ via base change along $B\hookrightarrow X$, we may
  assume that $X$ is a curve, and $\pi|_C:C\to X$ is a dominant morphism.  Consider
  the short exact sequence of sheaves
  \[
  0\lra \sha\lra\she\lra \shq\lra 0\ ,
  \]
  where $\shq\deq \she/\sha$ is a torsion sheaf on $X$ with support $T$.

  The vector bundle map $\pi^*\she\to\op$ is surjective, hence $\pi^*\sha$ maps
  surjectively onto a sub-line-bundle $\shb$ of $\op$. Since $\pi^*\sha$ is ample, so
  is the quotient $\shb$.
%}}
\end{proof}

\begin{comment}
  It is known that a direct sum of line bundles is L-big if and only if at least one
  of the summands is big; a direct sum of line bundles is V-big if all summands are
  big. What is the situation with AEAness?
\end{comment}

\sideremark{\sandor{060914}: This remark is a bit confusing. It needs a rewrite.}

\begin{remark}
  In the case of a line bundle $\shl$, the largest open subset over which the
  evaluation map
  \[
  \HH{0}{X}{\shl}\otimes\sO_X\lra \shl
  \]
  is surjective is the complement of the stable base locus $X\setminus
  \mathbb{B}(\shl)$.

  An L-big vector bundle $E$ is V-big if $\mathbb{B}_+(\op)$ is contained in a union
  of fibres over a proper Zariski closed subset of $X$.
\end{remark}

When $\shl$ is a line bundle on a surface, then it is immediate from the
intersection-theoretic characterizations that
\[
\mathbb  B_+(L) \equ \mathbb T(P_D)\ .
\]
Since we would need something along these lines on $\P(E)$, which in interesting cases has dimension at least three, the above observations can only serve as a pointer what kind of statements we would like to prove in higher dimensions.

\begin{lemma}\label{lemma-nef}
  Let $X$ be an irreducible projective variety, $\shl$ a nef line bundle on $X$. Then
  \[
  \mathbb B_+(\shl) \equ \mathbb T(\shl)\ .
  \]
\end{lemma}

\sideremark{\sandor{060914}: This proof needs a facelift!!}

\begin{proof}
  We have already observed in Remark~\ref{remark-augm} that $T(\shl)\subseteq
  B_+(\shl)$ in general.

  \begin{comment}
    According to the main result of \cite{Nak} (see also \cite[Corollary
    5.6]{ELMNP2}),
  \[
  B_+(\shl) \equ \bigcup_{V\subseteq X,\ (\shl^d\cdot V)=0}V\ ,
  \]
  where $V$ runs through all positive-dimensional irreducible subvarieties of $X$
  restricted to which $\shl$ ceases to be big; while
  \[
  \mathbb{T}(\shl) = T_{\eps} = \overline{\bigcup_{C\subseteq X,\ (\shl\cdot C)< \eps
      c_1(A)\cdot C } C } \textrm { for } 0< \eps \ll 1 ,
  \]
  where $C\subseteq X$ means an irreducible curve, and $A$ is a fixed ample divisor
  on $X$.  We have shown that $T(\shl)\subseteq B_+(\shl)$ in general.
  \end{comment}

To prove the other direction in the case where $L$ is nef,
we use  the main result of \cite{Nak} (see also \cite[Corollary 5.6]{ELMNP2}):
\[
 B_+(\shl) \equ \bigcup_{V\subseteq X,\ (\shl^d\cdot V)=0}V\ .
\]
We have to show that if
$V\subseteq X$ is an irreducible subvariety for which
  $\shl|_V$ is not big, then $V$ is contained in  $T(\shl)$. 
We will show that $V$ is covered by curves $C$ satisfying
  $\shl\cdot C < \eps A \cdot C$ for $\eps$ small enough.  

Assume first that $V$ is
  smooth, and apply \cite{BDPP} (the
  pseudo-effective cone is the dual of the cone of moving curves). 
Then
  $\shl|_V$ not big implies that it is not in the interior of 
the pseudo-effective  cone, 
hence there must exist a real 1-cycle $ 0 \neq C \in N_1(V)$ limit of moving
  curves $C_n\subseteq V$ with $\shl|_{V}\cdot C=0$.  
As $A|_{V}$ is ample on $V$,
  then $A|_{V} \cdot C >0$, 
so $\lim \frac{\shl|_{V}\cdot C_n}{A|_{V} \cdot C_n} =0$.
  Thus, for $n$ sufficiently large, $\shl \cdot C_n < \eps A \cdot C_n$.  The class
  of $C_n$ however covers $V$, which implies $V\subseteq \mathbb{T}(\shl)$.

  If $V$ is not smooth, then let $\mu:V'\to V \subseteq X$ be a resolution of
  singularities. Since $\shl|_V$ was pseudo-effective and not big to begin with, the
  same applies to $\mu^*(\shl|_V)$. Using the argument as above, there exist moving
  curves $C_n$ on $V'$ such that their limit is a non zero 1-cycle $C \in N_1(V')$
  such that $\mu^*(\shl|_V) \cdot C =0$, now let $A$ be an ample divisor on $X$, then
  $\mu^*(A)$ is big on $V'$, so $\mu^*(A) \cdot C >0$. Then as above $\lim
  \frac{\mu^*(\shl|_{V}) \cdot C_n}{\mu^*(A) \cdot C_n} =0$.  Thus, for $n$
  sufficiently large, $\mu^*(\shl|_{V}) \cdot C_n <\eps \mu^*(A) \cdot C_n$, hence by
  projection formula $\shl \cdot \mu_*(C_n) < \eps A \cdot \mu_*(C_n)$.  And the
  class of $\mu_*(C_n)$ covers $V$, which implies $V\subseteq \mathbb{T}(\shl)$.
\end{proof}

\begin{proposition}
  \label{aeavbig}
  Let $\she$ be a nef vector bundle on an irreducible projective variety $X$. Then
  $\she$ is AEA if and only if it is V-big.
\end{proposition}

\begin{proof}
  Both AEA and V-bigness imply that $\shl\deq \sO_\P(1)$ is a big and nef line
  bundle.  By the previous lemma, $B_+(\shl)=T(\shl)$, hence
  \[
  qE \text{ is V-big } \Leftrightarrow \pi(B_+(\shl))\neq X \Leftrightarrow
  \pi(T(\shl))\neq X \Leftrightarrow E \text{ is AEA}\ .\qedhere
  \]
\end{proof}

\begin{comment}
\silentremark{\sandor{060914}: Do we need this here?}
  Let now $D$ be a Cartier divisor on an irreducible projective variety $X$, and
  assume it has a Zariski decomposition; more precisely, I'd like to assume that
  there exists a unique minimal effective divisor $E$ on $X$ such that
  \[
  P \deq D-E
  \]
  is nef, and $P|_V$ is not big for any irreducible component $V$ of $E$. Following
  the argument of \cite[Example 1.11]{ELMNP}, we would get that
  \[
  B_+(\sO_X(D)) \equ \Null (P) \equ T(\sO_X(D))
  \]
  by Lemma~\ref{lemma-nef}.

\end{comment}

%*****************************************************************************

%
%
%\bibliographystyle{skalpha}
%\bibliography{bigne}

\begin{thebibliography}{ELMNP06}

\bibitem[BF12]{BauFun10}
{\sc T.~Bauer and M.~Funke}: \emph{Weyl and {Z}ariski chambers on {K}3
  surfaces}, Forum Math. \textbf{24} (2012), no.~3, 609--625. {\sf\scriptsize
  MR2926637}

\bibitem[BDPP13]{BDPP}
{\sc S.~Boucksom, J.-P. Demailly, M.~P{\u{a}}un, and T.~Peternell}: \emph{The
  pseudo-effective cone of a compact {K}\"ahler manifold and varieties of
  negative {K}odaira dimension}, J. Algebraic Geom. \textbf{22} (2013), no.~2,
  201--248. {\sf\scriptsize MR3019449}

\bibitem[CKL01]{cky}
{\sc D.~Cox, S.~Katz, and Y.-P.~Lee}, \emph{Virtual fundamental classes
  of zero loci}, Advances in algebraic geometry motivated by physics ({L}owell,
  {MA}, 2000), Contemp. Math., vol. 276, Amer. Math. Soc., Providence, RI,
  2001, pp.~157--166.

\bibitem[DR13]{divrous}
{\sc S.~{Diverio} and E.~{Rousseau}}: \emph{{The exceptional set and the
  Green-Griffiths locus do not always coincide}}, ArXiv e-prints (2013).

\bibitem[ELMNP06]{ELMNP}
{\sc L.~Ein, R.~Lazarsfeld, M.~Musta{\c{t}}{\u{a}}, M.~Nakamaye, and M.~Popa}:
  \emph{Asymptotic invariants of base loci}, Ann. Inst. Fourier (Grenoble)
  \textbf{56} (2006), no.~6, 1701--1734. {\sf\scriptsize MR2282673
  (2007m:14008)}

\bibitem[ELMNP09]{ELMNP2}
\bysame, \emph{Restricted volumes and base loci of linear series}, Amer. J.
  Math. \textbf{131} (2009), no.~3, 607--651.

\bibitem[Har70]{har70}
{\sc R.~Hartshorne}: \emph{Ample subvarieties of algebraic varieties}, Lecture
  Notes in Mathematics, Vol. 156, Springer-Verlag, Berlin-New York, 1970, Notes
  written in collaboration with C. Musili. {\sf\scriptsize MR0282977 (44
  \#211)}

\bibitem[Jab07]{KellyThesis}
{\sc K.~Jabbusch}: \emph{Notions of positivity for vector bundles}, Ph.D.
  thesis, University of Washington, 2007.

\bibitem[Lan86]{lang}
{\sc S.~Lang}: \emph{Hyperbolic and {D}iophantine analysis}, Bull. Amer. Math.
  Soc. (N.S.) \textbf{14} (1986), no.~2, 159--205. {\sf\scriptsize MR828820
  (87h:32051)}

\bibitem[Laz04]{PAG}
{\sc R.~Lazarsfeld}: \emph{Positivity in algebraic geometry. {I} , {II}},
  Ergebnisse der Mathematik und ihrer Grenzgebiete. 3. Folge. A Series of
  Modern Surveys in Mathematics [Results in Mathematics and Related Areas. 3rd
  Series. A Series of Modern Surveys in Mathematics], vol. 48, 49,
  Springer-Verlag, Berlin, 2004. {\sf\scriptsize MR2095471 (2005k:14001a)}

\bibitem[{Leh}11]{brian}
{\sc B.~{Lehmann}}: \emph{{The movable cone via intersections}}, ArXiv e-prints
  (2011).

\bibitem[{Les}12]{lesieutre}
{\sc J.~{Lesieutre}}: \emph{{The diminished base locus is not always closed}},
  ArXiv e-prints (2012).

\bibitem[Miy83]{Miy83}
{\sc Y.~Miyaoka}: \emph{Algebraic surfaces with positive indices},
  Classification of algebraic and analytic manifolds ({K}atata, 1982), Progr.
  Math., vol.~39, Birkh\"auser Boston, Boston, MA, 1983, pp.~281--301.
  {\sf\scriptsize MR728611 (85j:14067)}

\bibitem[Nak00]{Nak}
{\sc M.~Nakamaye}, \emph{Stable base loci of linear series}, Math. Ann.
  \textbf{318} (2000), no.~4, 837--847.

\bibitem[Nak04]{Nak04}
{\sc N.~Nakayama}: \emph{Zariski-decomposition and abundance}, MSJ Memoirs,
  vol.~14, Mathematical Society of Japan, Tokyo, 2004.

\bibitem[Urb07]{urbin}
{\sc S.~Urbinati}: \emph{On the properties of the integral part of ample and
  big divisors}, Master's thesis, Universit{\`a} degli Studi Roma Tre, 2007.

\bibitem[Vie83]{vie}
{\sc E.~Viehweg}: \emph{Weak positivity and the additivity of the {K}odaira
  dimension for certain fibre spaces}, Algebraic varieties and analytic
  varieties ({T}okyo, 1981), Adv. Stud. Pure Math., vol.~1, North-Holland,
  Amsterdam, 1983, pp.~329--353. {\sf\scriptsize MR715656 (85b:14041)}

\end{thebibliography}
%

%
%
%
%************************************************************************************************

%
%
%
%%

\end{document}